\documentclass{amsproc}

\def \Z{\Bbb Z}
\def \C{\Bbb C}

\def \N{\Bbb N}

\def \D{{\mathcal{D}}}
\def \E{{\mathcal{E}}}
\def \E{{\mathcal{E}}}

\def \S{{\mathcal{S}}}

\def \Res{{\rm Res}}

\def \End{{\rm End}}

\def \Hom{{\rm Hom}}

\def \<{\langle}
\def \>{\rangle}

\def \g{{\frak{g}}}

\def \be{\begin{equation}\label}
\def \ee{\end{equation}}
\def \bex{\begin{example}\label}
\def \eex{\end{example}}
\def \bl{\begin{lem}\label}
\def \el{\end{lem}}
\def \bt{\begin{thm}\label}
\def \et{\end{thm}}
\def \bp{\begin{prop}\label}
\def \ep{\end{prop}}
\def \br{\begin{rem}\label}
\def \er{\end{rem}}
\def \bc{\begin{coro}\label}
\def \ec{\end{coro}}
\def \bd{\begin{de}\label}
\def \ed{\end{de}}

\newtheorem{thm}{Theorem}[section]
\newtheorem{prop}[thm]{Proposition}
\newtheorem{coro}[thm]{Corollary}

\newtheorem{example}[thm]{Example}
\newtheorem{lem}[thm]{Lemma}
\newtheorem{rem}[thm]{Remark}
\newtheorem{de}[thm]{Definition}

\makeatletter \@addtoreset{equation}{section}

\makeatother \makeatletter

\begin{document}

\title{On quantum vertex algebras and their modules}

\author{Haisheng Li}
\address{Department of Mathematical Sciences, Rutgers University,
Camden, New Jersey 08102}
%\curraddr{Department of Mathematics and Statistics, Case Western
%Reserve University, Cleveland, Ohio 43403}
\email{hli@camden.rutgers.edu}
\thanks{The author was supported in part by NSF Grant DMS-0600189.}

%\author{Author Two}
%\address{Mathematical Research Section, School of Mathematical Sciences,
%Australian National University, Canberra ACT 2601, Australia}

%    General info
\subjclass{Primary 17B69, 17B68; Secondary 81R10}
%\date{January 1, 1994 and, in revised form, June 22, 1994.}

\dedicatory{To Geoffrey Mason on His 60th Birthday}

\keywords{Quantum vertex algebra, double Yangian,
Zamolodchikov-Faddeev algebra.}

\begin{abstract}
We give a survey on the developments in a certain theory of quantum
vertex algebras, including a conceptual construction of quantum
vertex algebras and their modules and a connection of double
Yangians and Zamolodchikov-Faddeev algebras with quantum vertex
algebras.
\end{abstract}

\maketitle

%\section*{This is an unnumbered first-level section head}
%This is an example of an unnumbered first-level heading.

%\specialsection*{This is a Special Section Head}
%This is an example of a special section head%
%%%%%%%%%%%%%%%%%%%%%%%%%%%%%%%%%%%%%%%%%%%%%%%%%%%%%%%%%%%%%%%%%%%%%%%%
%\footnote{Here is an example of a footnote. Notice that this footnote
%text is running on so that it can stand as an example of how a footnote
%with separate paragraphs should be written.
%\par
%And here is the beginning of the second paragraph.}%
%%%%%%%%%%%%%%%%%%%%%%%%%%%%%%%%%%%%%%%%%%%%%%%%%%%%%%%%%%%%%%%%%%%%%%%%
.
\section{Introduction}
In literature, there have been several notions of quantum vertex
algebra. The first, under the name ``deformed chiral algebra,'' was
introduced by E. Frenkel and N. Reshetikhin \cite{efr}, where
certain deformed Virasoro algebra and deformed $W$-algebras were
studied in this context. Then Etingof and Kazhdan \cite{ek}
introduced a notion of quantum vertex operator algebra in the sense
of formal deformation and they constructed a family of examples as
formal deformations of vertex operator algebras associated with
affine Lie algebras of type $A$. Later, Borcherds \cite{b-qva}
introduced another notion of quantum vertex algebra in terms of a
categorical language. These (three) notions, though closely related,
are quite different, and each of them leads to interesting algebraic
objects.

Inspired by the work of Etingof-Kazhdan, in a series of papers
(\cite{li-qva1}, \cite{li-qva2}, \cite{li-infinity}, \cite{kl}) we
have extensively studied a notion of (weak) quantum vertex algebra.
While quantum vertex operator algebras in the sense of \cite{ek} are
formal deformations of vertex algebras, (weak) quantum vertex
algebras are natural generalizations of vertex algebras and vertex
superalgebras. A weak quantum vertex algebra in this sense is a
vector space $V$ over $\C$ equipped with a linear map $Y:
V\rightarrow \Hom (V,V((x)))$ and equipped with a vector ${\bf 1}\in
V$, satisfying all the axioms that define the notion of a vertex
algebra except that the Jacobi identity axiom is replaced by an
$\S$-Jacobi identity axiom: For $u,v\in V$, there exists
$$\sum_{i=1}^{r}v^{(i)}\otimes u^{(i)}\otimes f_{i}(x)\in V\otimes
V\otimes \C((x))$$ such that
\begin{eqnarray*}
& &x_{0}^{-1}\delta\left(\frac{x_{1}-x_{2}}{x_{0}}\right)
Y(u,x_{1})Y(v,x_{2})\nonumber\\
&&\hspace{2cm}
-x_{0}^{-1}\delta\left(\frac{x_{2}-x_{1}}{-x_{0}}\right)
\sum_{i=1}^{r}f_{i}(-x_{0})Y(v^{(i)},x_{2})Y(u^{(i)},x_{1})\nonumber\\
& &=x_{2}^{-1}\delta\left(\frac{x_{1}-x_{0}}{x_{2}}\right)
Y(Y(u,x_{0})v,x_{2}).
\end{eqnarray*}
Just as Jacobi identity does,  $\S$-Jacobi identity implies
``associativity,'' so that weak quantum vertex algebras are nonlocal
vertex algebras, namely field algebras in the sense of \cite{bk},
which  are vertex analogs of noncommutative associative algebras.
The notion of quantum vertex algebra is a natural product of
Etingof-Kazhdna's notion of quantum vertex operator algebra, where a
quantum vertex algebra is a weak quantum vertex algebra $V$ equipped
with a unitary rational quantum Yang-Baxter operator $\S(x):
V\otimes V\rightarrow V\otimes V\otimes \C((x))$ such that for
$u,v\in V$, the element
$$\sum_{i=1}^{r}v^{(i)}\otimes u^{(i)}\otimes f_{i}(x)\in V\otimes
V\otimes \C((x))$$ for the $\S$-Jacobi identity is given by
$\S(x)(v\otimes u)$ and such that $\S(x)$ satisfies the
shift-condition and hexagon identity as in \cite{ek}. In this
theory, a notion of Etingof-Kazhdan, called ``non-degeneracy,''
plays a very important role just as it does in \cite{ek}. A result,
lifted from \cite{ek}, is that if a weak quantum vertex algebra $V$
is non-degenerate, then $V$ is a quantum vertex algebra with $\S(x)$
a linear map uniquely determined by the $\S$-Jacobi identity axiom.
In view of this, a natural way to construct quantum vertex algebras
is first to construct weak quantum vertex algebras and then to
establish non-degeneracy.

Non-degeneracy (in the sense of \cite{ek}) has been studied for
general vertex algebras in \cite{li-simple} and for general nonlocal
vertex algebras in \cite{li-qva2}. It was proved in \cite{li-simple}
that if $V$ is a vertex algebra of countable dimension (over $\C$)
and if $V$ as a $V$-module is irreducible, then $V$ is
non-degenerate.  In particular, this implies that every simple
vertex operator algebra in the sense of \cite{flm} is
non-degenerate. (Needless to say, this study was motivated by a
result of \cite{ek}. See Remark \ref{rnon-degenerate} for further
information.)  It was proved in \cite{li-qva2} that the same
assertion of \cite{li-simple} also holds for nonlocal vertex
algebras. Among other results on non-degeneracy, it was proved
therein that all the vertex algebras that are based on the canonical
universal vacuum modules for affine Lie algebras or for the Virasoro
Lie algebra are non-degenerate.

As for constructing weak quantum vertex algebras,  a conceptual
construction of weak quantum vertex algebras and their modules was
established in \cite{li-qva1}, generalizing the construction in
\cite{li-local} of vertex superalgebras and their modules (cf.
\cite{li-g1}, \cite{li-new}). Given a general vector space $W$, one
considers the space $\Hom (W,W((x)))$, alternatively denoted by
$\E(W)$. One then studies what were called compatible subsets and
$\S$-local subsets of $\E(W)$.  It was proved that every compatible
subset of $\E(W)$ generates a nonlocal vertex algebra with $W$ as
canonical (faithful) module, while every $\S$-local subset generates
a weak quantum vertex algebra with $W$ as a module. As an
application, examples of quantum vertex algebras were constructed by
using Zamolodchikov-Faddeev algebras (see \cite{li-qva2}, \cite{kl})
and certain versions of double Yangians (see \cite{li-infinity}).

Independently, Anguelova and Bergvelt \cite{ab} (cf. \cite{ang1},
\cite{ang2}) have studied what were called $H_{D}$-quantum vertex
algebras and they constructed very interesting examples. The work of
Anguelova and Bergvelt is closely related to both the work of
Borcherds and that of Etingof-Kazhdan. In particular, the notion of
$H_{D}$-quantum vertex algebra generalizes Etingof-Kazhdan's notion
of braided vertex operator algebra (see \cite{ek}) in a certain way.
One drawback of this generality is that $H_{D}$-quantum vertex
algebras {\em no longer} satisfy associativity, though they satisfy
a braided associativity (see \cite{ek}, \cite{ab}). Note that there
is a big difference in constructions; Anguelova and Bergvelt use
Borcherds' bicharacter-twisting method to establish $H_{D}$-quantum
vertex algebra structures on certain Hopf algebras, whereas we use
vertex operators on potential modules to construct quantum vertex
algebras, so that a representation theory is built up at the same
time.

A fundamental problem, which was posed by Igor Frenkel and Naihuan
Jing \cite{fj} in 1988 (cf. \cite{efk}), was to develop a theory of
quantum vertex algebras so that quantum affine algebras can be
associated with quantum vertex algebras in the same way that affine
Lie algebras were associated with vertex algebras. Arguably the
whole study on ``quantum vertex algebras'' was originated from this
problem. Though this particular problem is still to be solved, it
has led us to various interesting algebraic structures. In a recent
work \cite{li-hqva}, we have studied $\hbar$-adic quantum vertex
algebras and we have successfully associated centrally extended
double Yangians (see \cite{kh}) with $\hbar$-adic quantum vertex
algebras. It is our belief that one will be able to provide a
complete solution to the aforementioned problem, essentially in
terms of this theory of quantum vertex algebras.

We would like to thank Maarten Bergvelt, Gaywalee Yamskulna, and
Wenhua Zhao, for organizing this wonderful conference in honor of
Geoffrey Mason on the occasion of his 60th birthday, and we wish
Geoffrey a healthy and productive future for many years to come.

\section{Weak quantum vertex algebras and quantum vertex algebras}
In this section we shall discuss the notions of weak quantum vertex
algebra and quantum vertex algebra and the notion of nonlocal vertex
algebra. We shall also discuss Etingof-Kazhdan's notion of
non-degeneracy for nonlocal vertex algebras, which plays a crucial
role in the study of quantum vertex algebras. Throughout this paper,
$\N$ denotes the set of nonnegative integers and $\C$ denotes the
field of complex numbers, which will be the scalar field for vector
spaces.

We start with the notion of nonlocal vertex algebra (see
\cite{li-qva1}; cf. \cite{kac}, \cite{bk}, \cite{li-g1}).

\bd{dnonlocal-va} {\em A {\em nonlocal vertex algebra} is a vector
space $V$, equipped with a linear map
\begin{eqnarray*}
Y:&& V\rightarrow \Hom (V,V((x)))\subset (\End
V)[[x,x^{-1}]]\\
&&v\mapsto Y(v,x)=\sum_{n\in \Z}v_{n}x^{-n-1},
\end{eqnarray*} and equipped with
a vector ${\bf 1}\in V$, satisfying the conditions that $Y({\bf
1},x)=1$,
$$Y(v,x){\bf 1}\in V[[x]]\;\;\mbox{ and }\ \ (Y(v,x){\bf
1})|_{x=0}=v\ \ \mbox{ for }v\in V,$$ and that for any $u,v,w\in V$,
there exists $l\in \N$ such that
\begin{eqnarray}\label{ewassoc-def}
(x_{0}+x_{2})^{l}Y(u,x_{0}+x_{2})Y(v,x_{2})w =
(x_{0}+x_{2})^{l}Y(Y(u,x_{0})v,x_{2})w
\end{eqnarray}
(the {\em weak associativity}).} \ed

Note that the notion of nonlocal vertex algebra is the same as the
notion of weak $G_{1}$-vertex algebra in \cite{li-g1} and is
essentially the same as Bakalov and Kac's notion of field algebra
(see \cite{bk}; cf. \cite{kac}). It has been well understood that
vertex algebras are analogues and generalizations of commutative
associative algebras with identity. In contrast, nonlocal vertex
algebras are analogues and generalizations of (noncommutative)
associative algebras with identity.

For a nonlocal vertex algebra $V$, define a linear operator $\D$ on
$V$ by
\begin{eqnarray}
\D v=v_{-2}{\bf 1}=\frac{d}{dx}Y(v,x){\bf 1})|_{x=0}\ \ \ \mbox{ for
}v\in V.
\end{eqnarray}
Then
\begin{eqnarray}
[\D, Y(v,x)]=Y(\D v,x)=\frac{d}{dx}Y(v,x).
\end{eqnarray}

For a nonlocal vertex algebra $V$, the vertex operator map $Y:
V\rightarrow \Hom (V,V((x)))$ amounts to a linear map from $V\otimes
V$ to $V((x))$, which we denote by
$$Y(x): V\otimes V \rightarrow V((x)),$$
following \cite{ek}. The following notion of weak quantum vertex
algebra, which was formulated in \cite{li-qva1}, singles out a
special family of nonlocal vertex algebras:

\bd{dwqva} {\em  A {\em weak quantum vertex algebra} is a nonlocal
vertex algebra $V$ that satisfies {\em $\S$-locality}: For $u,v\in
V$, there exists
$$\sum_{i=1}^{r}v^{(i)}\otimes u^{(i)}\otimes f_{i}(x)\in V\otimes
V\otimes \C((x))$$ such that
\begin{eqnarray}\label{es-locality-def}
&&(x_{1}-x_{2})^{k}Y(u,x_{1})Y(v,x_{2})=(x_{1}-x_{2})^{k}
\sum_{i=1}^{r}f_{i}(x_{2}-x_{1})Y(v^{(i)},x_{2})Y(u^{(i)},x_{1})\nonumber\\
&&
\end{eqnarray}
for some nonnegative integer $k$, where $f_{i}(x_{2}-x_{1})$ is to
be expanded in the nonnegative powers of $x_{1}$, that is
$$f_{i}(x_{2}-x_{1})=e^{-x_{1}\frac{\partial}{\partial
x_{2}}}f_{i}(x_{2}),$$ in view of the formal Taylor theorem.} \ed

\br{rwqva-wc-wa} {\em Notice that if $(V,Y,{\bf 1})$ satisfies all
the axioms that define the notion of nonlocal vertex algebra except
the weak associativity axiom, then  for $u,v,w\in V$ the
$\S$-locality relation (\ref{es-locality-def}) together with the
weak associativity relation (\ref{ewassoc-def}) is equivalent to
\begin{eqnarray}\label{esjacobi-def}
& &x_{0}^{-1}\delta\left(\frac{x_{1}-x_{2}}{x_{0}}\right)
Y(u,x_{1})Y(v,x_{2})w\nonumber\\
&&\hspace{2cm}
-x_{0}^{-1}\delta\left(\frac{x_{2}-x_{1}}{-x_{0}}\right)
\sum_{i=1}^{r}f_{i}(-x_{0})Y(v^{(i)},x_{2})Y(u^{(i)},x_{1})w\nonumber\\
& &=x_{2}^{-1}\delta\left(\frac{x_{1}-x_{0}}{x_{2}}\right)
Y(Y(u,x_{0})v,x_{2})w
\end{eqnarray}
(the {\em $\S$-Jacobi identity}). In view of this, the notion of
weak quantum vertex algebra can be defined by using the axioms which
define the notion of nonlocal vertex algebra with weak associativity
replaced by $\S$-Jacobi identity.} \er

The notion of weak quantum vertex algebra clearly generalizes the
notion of vertex algebra and the notion of vertex superalgebra. A
quantum vertex algebra will be a weak quantum vertex algebra for
which the $\S$-locality is controlled by a unitary rational quantum
Yang-Baxter operator as we recall next.

\bd{dqybo} {\em Let $U$ be a vector space. A {\em unitary rational
quantum Yang-Baxter operator with one spectral parameter} on $U$ is
a linear map
$$\S(x): U\otimes U\rightarrow U\otimes U\otimes \C((x))$$
such that
\begin{eqnarray}
&&\S(-x)\S_{21}(x)=1,\\
&&\S_{12}(x)\S_{13}(x+z)\S_{23}(z)=\S_{23}(z)\S_{13}(x+z)\S_{12}(x),
\end{eqnarray}
where $\S_{12}(x), \S_{13}(x), \S_{23}(x)$ are linear maps from
$U\otimes U\otimes U$ to $U\otimes U\otimes U\otimes \C((x))$, e.g.,
$$\S_{12}(x)=\S(x)\otimes 1,\ \ \S_{23}(x)=1\otimes \S(x),$$
and where $\S_{21}(x)=P\S_{12}(x)P$ with $P$ the flip on $U\otimes
U$, i.e., $P(u\otimes v)=v\otimes u$ for $u,v\in U$. } \ed

The notion of quantum vertex algebra (cf. \cite{li-qva1}) was
derived from the notion of quantum vertex operator algebra in
\cite{ek}.

\bd{dqva} {\em A {\em quantum vertex algebra} is a weak quantum
vertex algebra $V$ equipped with a unitary rational quantum
Yang-Baxter operator $\S(x): V\otimes V\rightarrow V\otimes V\otimes
\C((x))$ such that for $u,v\in V$, the ${\mathcal{S}}$-Jacobi
identity (\ref{esjacobi-def}) holds with
$$\sum_{i=1}^{r}v^{(i)}\otimes u^{(i)}\otimes f_{i}(x)
={\mathcal{S}}(x)(v\otimes u)$$ and such that
\begin{eqnarray}
& &[\D\otimes 1,\S(x)]=-\frac{d}{dx}\S(x),\\
& &\S(z)(Y(x)\otimes 1) =(Y(x)\otimes 1)\S_{23}(z)\S_{13}(z+x),
\end{eqnarray}
called the {\em shift-condition} and {\em hexagon identity},
respectively.} \ed

In their study on quantum vertex operator algebras, Etingof and
Kazhdan introduced in \cite{ek} a very important notion, called
non-degeneracy. A weak quantum vertex algebra (or more generally a
nonlocal vertex algebra) $V$ is said to be {\em non-degenerate} if
for every positive integer $n$, the linear map
\begin{eqnarray*}
Z_{n}: V^{\otimes n}\otimes \C((x_{1}))\cdots ((x_{n})) \rightarrow
V((x_{1}))\cdots ((x_{n})),
\end{eqnarray*}
defined by
\begin{eqnarray*}
Z_{n}(v^{(1)}\otimes \cdots \otimes v^{(n)}\otimes
f(x_{1},\dots,x_{2}))=f(x_{1},\dots,x_{n})Y(v^{(1)},x_{1})\cdots
Y(v^{(n)},x_{n}){\bf 1},
\end{eqnarray*}
is injective. The importance of this notion can be seen from the
following result, which, proved in \cite{li-qva1}, is lifted from
Proposition 1.11 of \cite{ek}:

\bp{pek} Let $V$ be a weak quantum vertex algebra. If $V$ is
non-degenerate, then there exists a unique linear map $\S(x):
V\otimes V\rightarrow V\otimes V\otimes \C((x))$ such that for
$u,v\in V$, the ${\mathcal{S}}$-Jacobi identity (\ref{esjacobi-def})
holds with
$$\sum_{i=1}^{r}v^{(i)}\otimes u^{(i)}\otimes f_{i}(x)
={\mathcal{S}}(x)(v\otimes u).$$ Furthermore, $(V,\S(x))$ is a
quantum vertex algebra. \ep

\br{rnon-degenerate} {\em Let $\g$ be a finite-dimensional
semisimple Lie algebra and let $\hat{\g}$ be the corresponding
(untwisted) affine Kac-Moody Lie algebra. For every complex number
$\ell$, one has a vertex algebra $V_{\hat{\g}}(\ell,0)$, based on
the canonical universal vacuum $\hat{\g}$-module of level $\ell$.
 It was proved in \cite{ek} that if $V_{\hat{\g}}(\ell,0)$
as a $\hat{\g}$-module is irreducible, then $V_{\hat{\g}}(\ell,0)$
is non-degenerate. With this as the main motivation, it was proved
in \cite{li-simple} that if a vertex algebra $V$ is of countable
dimension (over $\C$) and if $V$ as a $V$-module is irreducible,
then $V$ is non-degenerate. This in particular implies that every
simple vertex operator algebra in the sense of \cite{flm} is
non-degenerate. Furthermore, it was proved in \cite{li-qva2} that if
a nonlocal vertex algebra $V$ is of countable dimension (over $\C$)
and if $V$ as a $V$-module is irreducible, then $V$ is
non-degenerate.} \er

As every non-degenerate weak quantum vertex algebra has a unique
quantum vertex algebra structure by Proposition \ref{pek}, the term
{\em ``non-degenerate quantum vertex algebra''} is unambiguous.
Similarly, the term {\em ``irreducible quantum vertex algebra''} is
meaningful.

Having established non-degeneracy for irreducible nonlocal vertex
algebras, we consider non-degeneracy for reducible nonlocal vertex
algebras. A simple fact is that every nonlocal vertex subalgebra of
a non-degenerate nonlocal vertex algebra is always non-degenerate.
In view of this, a nonlocal vertex algebra $V$ may be still
non-degenerate even though $V$ as a $V$-module is reducible. On the
other hand, this gives us a practical way to establish
non-degeneracy
--- by embedding a reducible nonlocal vertex algebra into an
irreducible nonlocal vertex algebra.

There is a simple construction of (nonlocal) vertex algebras, due to
Borcherds (see \cite{b-va}), based on classical associative algebras
with derivations. Let $A$ be an associative algebra with identity
$1$ equipped with a derivation $d$. Define a linear map $Y:
A\rightarrow (\End_{\C}A)[[x]]$ by
$$Y(a,x)b=(e^{xd}a)b\;\;\;\mbox{ for }a,b\in A.$$
Then $(A,Y,1)$ is a nonlocal vertex algebra (with $1$ as the vacuum
vector). Furthermore, this nonlocal vertex algebra is a vertex
algebra if and only if $A$ is commutative. Denote this nonlocal
vertex algebra by $V(A,d)$. The fact is that any nonlocal vertex
algebra $V$, satisfying the condition that $Y(u,x)v\in V[[x]]$ for
$u,v\in V$, is isomorphic to $V(A,d)$ with $A=V$, $u\cdot
v=(Y(u,x)v)|_{x=0}=u_{-1}v$ for $u,v\in V$, and $d=\D$. We shall see
that this construction gives us many important reducible but
non-degenerate vertex algebras.

\br{rfree-nondeg} {\em  Let $U$ be a finite-dimensional vector
space. Consider the symmetric algebra $S(U\otimes
t^{-1}\C[t^{-1}])$. Equipping this commutative associative algebra
with a derivation $d/dt$, we obtain a vertex algebra $V(S(U\otimes
t^{-1}\C[t^{-1}]),d/dt)$. By embedding this vertex algebra  into a
simple vertex operator algebra (associated to a Heisenberg Lie
algebra), it was proved in \cite{li-qva2} that $V(S(U\otimes
t^{-1}\C[t^{-1}]),d/dt)$ is non-degenerate.} \er

\br{rnon-degenerate-filtration} {\em Let $V$ be a general nonlocal
vertex algebra. Let $\E=\{E_{n}\}_{n\in \Z}$ be an increasing
filtration of subspaces of $V$, satisfying the condition that ${\bf
1}\in E_{0}$,
$$Y(u,x)E_{n}\subset E_{m+n}[[x,x^{-1}]]\ \ \
\mbox{ for }m,n\in \Z,\; u\in E_{m}.$$ Set $Gr_{\E}(V)=\oplus_{n\in
\Z} (E_{n}/E_{n-1})$, a vector space. Then (cf. \cite{li-vpa})
$Gr_{\E}(V)$ is a nonlocal vertex algebra with ${\bf 1}+E_{-1}$ as
the vacuum vector and with
$$Y(u+E_{m-1},x)(v+E_{n-1})=\sum_{k\in \Z}(u_{k}v+E_{m+n-1})x^{-k-1}\ \ \
\mbox{ for }u\in E_{m},\; v\in E_{n}.$$ It was proved in
\cite{li-qva2} that if the filtration $\E=\{E_{n}\}_{n\in \Z}$ is
lower truncated in the sense that $E_{n}=0$ for $n$ sufficiently
negative and if $Gr_{\E}(V)$ is non-degenerate, then $V$ is
non-degenerate. By using this result it was proved therein that the
vertex algebras associated to affine Lie algebras and the Virasoro
algebra, based on the canonical universal vacuum modules (cf. Remark
\ref{rnon-degenerate}), are non-degenerate.} \er

\bex{example-lie} {\em  Let $\g$ be any finite-dimensional Lie
algebra over $\C$. Consider the Lie algebra $t^{-1}\g[t^{-1}]$
$(=\g\otimes t^{-1}\C[t^{-1}])$. As $d/dt$ acts on
$t^{-1}\g[t^{-1}]$ as a derivation, it canonically lifts to a
derivation of the universal enveloping algebra
$U(t^{-1}\g[t^{-1}])$. Then we have a nonlocal vertex algebra
$V(U(t^{-1}\g[t^{-1}]), d/dt)$. For $a\in \g$, we have
$$Y(a\otimes t^{-1},x)=e^{x(d/dt)}(a\otimes t^{-1})=a\otimes
(t+x)^{-1}=\sum_{n\ge 0}(a\otimes t^{-n-1})x^{n}.$$ Denote
$\sum_{n\ge 0}(a\otimes t^{-n-1})x^{n}$ by $a(x)^{-}$, which is the
singular part of the whole field $$a(x)=\sum_{n\in \Z}(a\otimes
t^{n})x^{-n-1}.$$ For $a,b\in \g$, we have (cf. \cite{efk})
$$[a(x_{1})^{-},b(x_{2})^{-}]
=(x_{2}-x_{1})^{-1}([a,b](x_{1})^{-}-[a,b](x_{2})^{-}),$$ or
equivalently,
\begin{eqnarray}\label{eab-singular}
\ \ \ \ \ \ \ \ \ a(x_{1})^{-}b(x_{2})^{-}=b(x_{2})^{-}a(x_{1})^{-}
+(x_{2}-x_{1})^{-1}([a,b](x_{1})^{-}-[a,b](x_{2})^{-}).
\end{eqnarray}
 Using this and induction one can show that
$V(U(t^{-1}\g[t^{-1}]), d/dt)$ is a weak quantum vertex algebra.  By
using Remarks \ref{rfree-nondeg} and
\ref{rnon-degenerate-filtration} it was proved in \cite{li-qva2}
that $V(U(t^{-1}\g[t^{-1}]),d/dt)$, which is a genuine nonlocal
vertex algebra, is in fact a non-degenerate quantum vertex algebra
with $\g\otimes t^{-1}$ as a generating subspace. Note that the
whole fields $a(x)$ for $a\in \g$ satisfy the locality relation
$$(x_{1}-x_{2})^{2}a(x_{1})b(x_{2})=(x_{1}-x_{2})^{2}b(x_{2})a(x_{1})$$
(with trivial braiding) and they generate the vertex algebra
$V_{\hat{\g}}(\ell,0)$ as in Remark \ref{rnon-degenerate}, where the
half fields $a(x)^{-}$ for $a\in \g$ satisfy $\S$-locality relation
(\ref{eab-singular}) with nontrivial braiding and they generate the
quantum vertex algebra $V(U(t^{-1}\g[t^{-1}]), d/dt)$.} \eex

\section{General construction of weak quantum vertex algebras}
We have seen from Proposition \ref{pek} that one way to construct a
quantum vertex algebra is first to construct a weak quantum vertex
algebra and then to establish non-degeneracy. Here we discuss a
conceptual construction of weak quantum vertex algebras and their
modules.

Let us begin with the notions of module and quasi module for a
nonlocal vertex algebra (see \cite{li-qva1}).

\bd{dmodule} {\em Let $V$ be a general nonlocal vertex algebra. A
{\em $V$-module} is a vector space $W$, equipped with a linear map
$$Y_{W}: V\rightarrow \Hom (W,W((x)))\subset (\End W)[[x,x^{-1}]],$$
satisfying the conditions that $Y_{W}({\bf 1},x)=1_{W}$ (the
identity operator on $W$) and that for $u,v\in V,\; w\in W$, there
exists a nonnegative integer $l$ such that
$$(x_{0}+x_{2})^{l}Y_{W}(u,x_{0}+x_{2})Y_{W}(v,x_{2})w
=(x_{0}+x_{2})^{l}Y_{W}(Y(u,x_{0})v,x_{2})w.$$ In this definition,
if we replace the second condition with ``for $u,v\in V,\; w\in W$,
there exists a nonzero polynomial $p(x_{1},x_{2})$ such that
$$p(x_{0}+x_{2},x_{2})Y_{W}(u,x_{0}+x_{2})Y_{W}(v,x_{2})w
=p(x_{0}+x_{2},x_{2})Y_{W}(Y(u,x_{0})v,x_{2})w,\mbox{''}$$ we then
arrive at the notion of {\em quasi $V$-module}.} \ed

A conceptual construction of weak quantum vertex algebras and their
modules, which was established in \cite{li-qva1}, is to use formal
vertex operators of a certain type. Let $W$ be a general vector
space over $\C$. Set
$$\E(W)=\Hom (W,W((x))),$$
which can be viewed as a vertex analog of the notion of $\End W$.
The identity operator, which is denoted by $1_{W}$, is a typical
element of $\E(W)$. An (ordered) sequence
$(\psi^{(1)}(x),\dots,\psi^{(r)}(x))$ in $\E(W)$ is said to be {\em
quasi-compatible} if there exists a nonzero polynomial
$p(x_{1},x_{2})$ such that
\begin{eqnarray}
\ \ \ \ \ \ \ \ \left(\prod_{1\le i<j\le r}p(x_{i},x_{j})\right)
\psi^{(1)}(x_{1})\cdots \psi^{(r)}(x_{r})\in \Hom
(W,W((x_{1},\dots,x_{r}))).
\end{eqnarray}
A subset $U$ of $\E(W)$ is said to be {\em quasi-compatible} if
every finite sequence in $U$ is quasi-compatible. In this
definition, if we assume that $p(x_{1},x_{2})$ is of the form
$(x_{1}-x_{2})^{k}$ with $k\in \N$, we then arrive at the notion of
{\em compatibility}.

Let $\C(x)$ denote the field of rational functions. Let
\begin{eqnarray}
\iota_{x;0}:\C(x)\rightarrow \C((x))
\end{eqnarray}
be the unique field embedding that preserves each polynomial. Given
any nonzero polynomial $q(x)$, writing $q(x)=x^{k}(\alpha +xh(x))$
for some $k\in \N,\; \alpha \in \C,\;h(x)\in \C[x]$ with $\alpha\ne
0$, we have
\begin{eqnarray}
\iota_{x;0}\left(\frac{1}{q(x)}\right)=\sum_{j\ge
0}\alpha^{-j-1}x^{j-k}h(x)^{j}\in \C((x)).
\end{eqnarray}
Furthermore, let $\C(x_{1},x_{2})$ denote the field of rational
functions. As $\C[x_{1},x_{2}]$ is a subring of the field
$\C((x_{1}))((x_{2}))$, there exists a unique field-embedding
\begin{eqnarray}
\iota_{x_{1},x_{2}}:\ \
 \C(x_{1},x_{2})\rightarrow \C((x_{1}))((x_{2})),
\end{eqnarray}
preserving each polynomial. For any nonzero polynomial
$p(x_{1},x_{2})$, writing
$p(x_{1},x_{2})=x_{2}^{k}(f(x_{1})+x_{2}g(x_{1},x_{2}))$ for some
$k\in \N,\; f(x)\in \C[x],\; g(x_{1},x_{2})\in \C[x_{1},x_{2}]$ with
$f(x)\ne 0$, we have
\begin{eqnarray}
\iota_{x_{1},x_{2}}\left(\frac{1}{p(x_{1},x_{2})}\right)=\sum_{j\ge
0}\iota_{x_{1};0}(f(x_{1})^{-j-1})x_{2}^{j-k}g(x_{1},x_{2})^{j}.
\end{eqnarray}

Let $a(x),b(x)\in \E(W)$ such that the ordered pair $(a(x),b(x))$ is
quasi compatible. By definition, there exists $0\ne
p(x_{1},x_{2})\in \C[x_{1},x_{2}]$ such that
\begin{eqnarray}\label{epab-condition}
p(x_{1},x_{2})a(x_{1})b(x_{2})\in \Hom (W,W((x_{1},x_{2}))).
\end{eqnarray}
Note that for any formal series $A(x_{1},x_{2})\in \Hom
(W,W((x_{1},x_{2})))$, $A(x_{2},x_{2})$ exists in $\Hom
(W,W((x_{2})))$. Furthermore, $A(x_{2}+x_{0},x_{2})$, which is equal
to
$$\left(e^{x_{0}\frac{\partial}{\partial
x_{1}}}A(x_{1},x_{2})\right)|_{x_{1}=x_{2}}$$ lies in $(\Hom
(W,W((x_{2}))))[[x_{0}]]$. In view of this, the expression
$$\iota_{x,x_{0}}\left(\frac{1}{p(x_{0}+x,x)}\right)
\left(p(x_{1},x)a(x_{1})b(x)\right)|_{x_{1}=x+x_{0}}$$  exists in
$(\Hom (W,W((x_{2}))))[[x_{0}]]$.

\bd{danb} {\em Let $a(x),b(x)\in \E(W)$ such that $(a(x),b(x))$ is
quasi compatible. Define $a(x)_{n}b(x)\in \E(W)$ for $n\in \Z$ in
terms of generating function
\begin{eqnarray*}
Y_{\E}(a(x),x_{0})b(x) =\sum_{n\in \Z}a(x)_{n}b(x) x_{0}^{-n-1}
\end{eqnarray*}
by
\begin{eqnarray}
\ \ \ \ \ \ \
Y_{\E}(a(x),x_{0})b(x)=\iota_{x,x_{0}}\left(\frac{1}{p(x_{0}+x,x)}\right)
\left(p(x_{1},x)a(x_{1})b(x)\right)|_{x_{1}=x+x_{0}},
\end{eqnarray}
where $p(x_{1},x_{2})$ is any nonzero polynomial such that
(\ref{epab-condition}) holds.} \ed

This defines (binary) partial operations on $\E(W)$ parameterized by
$n\in \Z$. A quasi-compatible subspace $U$ of $\E(W)$ is said to be
{\em $Y_{\E}$-closed} if
\begin{eqnarray*}
a(x)_{n}b(x)\in U\;\;\;\mbox{ for }a(x),b(x)\in U,\; n\in \Z.
\end{eqnarray*}

We have the following conceptual results (see \cite{li-qva1}; cf.
\cite{li-g1}):

\bt{tqmodule} Let $W$ be a general vector space. Assume that $V$ is
a $Y_{\E}$-closed quasi-compatible subspace of $\E(W)$, containing
$1_{W}$. Then $(V,Y_{\E},1_{W})$ carries the structure of a nonlocal
vertex algebra with $W$ as a canonical quasi-module where
$$Y_{W}(\alpha(x),x_{0})=\alpha(x_{0})\ \ \mbox{ for }\alpha(x)\in V.$$
Furthermore, if $V$ is compatible, then $W$ is a $V$-module. \et

\bt{tmaximal} Let $W$ be a general vector space and let $U$ be a
quasi-compatible subset of $\E(W)$. Then there exists a (unique)
smallest $Y_{\E}$-closed subspace $\<U\>$ of $\E(W)$, containing $U$
and $1_{W}$. Furthermore,  $(\<U\>,Y_{\E},1_{W})$ carries the
structure of a nonlocal vertex algebra with $W$ as a quasi-module
where $Y_{W}(\alpha(x),x_{0})=\alpha(x_{0})$ for $\alpha(x)\in
\<U\>$. \et

To obtain a weak quantum vertex algebra, one needs to consider
subsets of a certain type of $\E(W)$.

\bd{ds-locality}{\em A subspace $U$ of $\E(W)$ is said to be {\em
${\mathcal{S}}$-local} if for any $a(x),b(x)\in U$, there exists
$$\sum_{i=1}^{r}b^{(i)}(x)\otimes a^{(i)}(x)\otimes f_{i}(x)\in
U\otimes U\otimes \C((x))$$ such that
\begin{eqnarray}\label{es-local-def}
\ \ \ \ (x_{1}-x_{2})^{k}a(x_{1})b(x_{2})=(x_{1}-x_{2})^{k}
\sum_{i=1}^{r}f_{i}(x_{2}-x_{1})b^{(i)}(x_{2})a^{(i)}(x_{1})
\end{eqnarray}
for some nonnegative integer $k$. Furthermore, a subset is said to
be {\em $\S$-local} if its linear span is $\S$-local.} \ed

\br{rdifferent} {\em Let $a(x),b(x)\in \E(W)$. Assume that there
exists
$$\sum_{i=1}^{r}b^{(i)}(x)\otimes a^{(i)}(x)\otimes f_{i}(x)\in
\E(W)\otimes \E(W)\otimes \C((x))$$ such that (\ref{es-local-def})
holds. One can show (see \cite{li-qva1}) that $(a(x),b(x))$ is
compatible and
\begin{eqnarray*}
&&Y_{\E}(a(x),x_{0})b(x)\nonumber\\
&=&\Res_{x_{1}}x_{0}^{-1}\delta\left(\frac{x_{1}-x}{x_{0}}\right)a(x_{1})b(x)
-x_{0}^{-1}\delta\left(\frac{x-x_{1}}{-x_{0}}\right)
\sum_{i=1}^{r}f_{i}(x-x_{1})b^{(i)}(x)a^{(i)}(x_{1}).
\end{eqnarray*}
In terms of components, for $n\in \Z$ we have
\begin{eqnarray*}
&&a(x)_{n}b(x)\nonumber\\
&=&\Res_{x_{1}}\left((x_{1}-x)^{n}a(x_{1})b(x)
-\sum_{i=1}^{r}(-x+x_{1})^{n}f_{i}(x-x_{1})b^{(i)}(x)a^{(i)}(x_{1})\right).
\end{eqnarray*}}
\er

We have the following general result (see \cite{li-qva1}):

\bt{ts-local} Let $U$ be an $\S$-local subset of $\E(W)$. Then $U$
is compatible and $\<U\>$ is the smallest $Y_{\E}$-closed $\S$-local
subspace of $\E(W)$, containing $U$ and $1_{W}$. Furthermore,
$(\<U\>,Y_{\E},1_{W})$ is a weak quantum vertex algebra with $W$ as
a module. \et

The following result (see \cite{li-qva1}) should be used as a
companion of Theorem \ref{ts-local} for practical applications (cf.
Remark \ref{rtransport}):

\bp{ptransported} Let $W$ be a general vector space, let $V$ be a
weak quantum vertex algebra generated by an $\S$-local subset of
$\E(W)$, and let
$$u(x),v(x),u^{(1)}(x),v^{(1)}(x),\dots,u^{(r)}(x),v^{(r)}(x),c^{0}(x),\dots,c^{s}(x) \in V,$$
$$f_{1}(x),\dots,f_{r}(x)\in \C((x)).$$ Suppose that
\begin{eqnarray}\label{e-cross-bracket-W}
& &(x_{1}-x_{2})^{n}u(x_{1})v(x_{2})
-(-x_{2}+x_{1})^{n}\sum_{i=1}^{r}f_{i}(x_{2}-x_{1})v^{(i)}(x_{2})u^{(i)}(x_{1})
\nonumber\\
&&\ \ \ \ =
\sum_{j=0}^{s}c^{j}(x_{2})\frac{1}{j!}\left(\frac{\partial}{\partial
x_{2}}\right)^{j} x_{2}^{-1}\delta\left(\frac{x_{1}}{x_{2}}\right)
\end{eqnarray}
for some integer $n$, then
\begin{eqnarray}\label{e-cross-bracket-algebra-case}
& &(x_{1}-x_{2})^{n}Y_{\E}(u(x),x_{1})Y_{\E}(v(x),x_{2})\nonumber\\
& &\ \ \ \ \ \ \ \ \
-(-x_{2}+x_{1})^{n}\sum_{i=1}^{r}f_{i}(x_{2}-x_{1})
Y_{\E}(v^{(i)}(x),x_{2})Y_{\E}(u^{(i)}(x),x_{1})
\nonumber\\
&&\ \ \ \
=\sum_{j=0}^{s}Y_{\E}(c^{j}(x),x_{2})\frac{1}{j!}\left(\frac{\partial}{\partial
x_{2}}\right)^{j} x_{2}^{-1}\delta\left(\frac{x_{1}}{x_{2}}\right).
\end{eqnarray}
\ep

\br{rtransport} {\em Here we outline a general procedure to
construct weak quantum vertex algebras. In reality, one often starts
with a certain associative algebra say $A$, which is represented by
generators and defining relations, such that generating functions
for generators form an $\S$-local set say $U$. For any highest
weight type $A$-module $W$, $U$ becomes an $\S$-local subset of
$\E(W)$, thus by Theorem \ref{ts-local} $U$ generates a weak quantum
vertex algebra $\<U\>_{W}$ with $W$ as a canonical module. Then
using Proposition \ref{ptransported} one shows that $\<U\>_{W}$ is a
vacuum $A$-module with each generating function $a(z)$ acting as
$Y_{\E}(a(x),z)$. Now, taking $W=V$ to be a universal ``vacuum
$A$-module,'' one shows that $\<U\>_{V}$ as a vacuum $A$-module is
isomorphic to $V$. Consequently, $V$ has a weak quantum vertex
algebra structure transported from $\<U\>_{V}$. Furthermore, for a
general highest weight type $A$-module $W$, since $\<U\>_{W}$ is a
vacuum $A$-module and since $V$ is a universal vacuum $A$-module,
there exists an $A$-module homomorphism $\psi_{W}$ from $V$ onto
$\<U\>_{W}$, sending the vacuum vector of $V$ to $1_{W}$. It follows
that $\psi_{W}$ is in fact a homomorphism of weak quantum vertex
algebras. As $W$ is a canonical module for $\<U\>_{W}$, $W$ is
naturally a $V$-module. } \er

\br{rbraiding-calculation} {\em As it was illustrated in Remark
\ref{rtransport}, we often obtain a non-degenerate quantum vertex
algebra structure on a space $V$ with an $\S$-local set of fields on
$V$, parameterized by the vectors of a subset $U$ of $V$. In this
case, $U$ becomes a generating subset of $V$ as a quantum vertex
algebra. The braiding for the generators is simply coded in the
$\S$-locality relation, while the braiding for the descendants can
be obtained by using the hexagon identity and unitarity relation. }
\er

\section{Examples of quantum vertex algebras}
Here we discuss some examples of quantum vertex algebras, which were
constructed by using the procedure outlined in Remark
\ref{rtransport}. The associative algebras that we are going to use
are Zamolodchikov-Faddeev algebras (see \cite{zam}, \cite{fad}) and
(centerless) double Yangians (see \cite{kt}).

Let $H$ be a finite-dimensional vector space over $\C$, equipped
with a bilinear form $\<\cdot,\cdot\>$, and let
$$\S(x): H\otimes H\rightarrow H\otimes H\otimes \C((x))$$
be a linear map. Let $T(H\otimes \C[t,t^{-1}])$ denote the free
associative algebra over $H\otimes \C[t,t^{-1}]$. For $u\in H,\;
n\in \Z$, set $u(n)=u\otimes t^{n}$. Form the generating function
$$u(x)=\sum_{n\in \Z}u(n)x^{-n-1}\in T(H\otimes
\C[t,t^{-1}])[[x,x^{-1}]].$$ Define an {\em $(H,\S)$-module} to be a
$T(H\otimes \C[t,t^{-1}])$-module $W$ satisfying the conditions that
for every $u\in H,\; w\in W$, $u(n)w=0$ for $n$ sufficiently large
and that for $u,v\in H,\; w\in W$,
\begin{eqnarray}\label{e2.18}
\ \ \ \ \ \ \ \
u(x_{1})v(x_{2})w-\sum_{i=1}^{r}f_{i}(x_{2}-x_{1})v^{(i)}(x_{2})u^{(i)}(x_{1})w
=\<u,v\>x_{2}^{-1}\delta\left(\frac{x_{1}}{x_{2}}\right)w,
\end{eqnarray}
where $\S(x)(v\otimes u)=\sum_{i=1}^{r}v^{(i)}\otimes u^{(i)}\otimes
f_{i}(x)$.  A {\em vacuum $(H,\S)$-module} is an $(H,\S)$-module $W$
equipped with a vector $w_{0}$ satisfying the conditions that
$u(n)w_{0}=0$ for $u\in H,\; n\ge 0$ and that $W=T(H\otimes
\C[t,t^{-1}])w_{0}$. In terms of this notion we have (see
\cite{li-qva2}):

\bt{tzf} Let $(V,{\bf 1})$ be a  vacuum $(H,\S)$-module that is
universal in the obvious sense. There exists a weak quantum vertex
algebra structure on $V$ with ${\bf 1}$ as the vacuum vector and
with
$$Y(u(-1){\bf
1},x)=u(x)\ \ \ \mbox{ for }u\in H.$$ Furthermore, for any
$(H,\S)$-module $W$, there exists a $V$-module structure such that
$Y_{W}(u(-1){\bf 1},x)=u(x)$ for $u\in H$. \et

Let $J_{+}$ denote the left ideal of $T(H\otimes \C[t,t^{-1}])$,
generated by the elements
$$u^{(1)}(n_{1})\cdots u^{(r)}(n_{r})$$
for $r\ge 1,\; u^{(i)}\in H,\; n_{i}\in \Z$ with $n_{1}+\cdots
+n_{r}\ge 0$. Set $$\widetilde{V}(H,\S)=T(H\otimes
\C[t,t^{-1}])/J_{+},$$ a $T(H\otimes \C[t,t^{-1}])$-module. For any
$u\in H,\; w\in \widetilde{V}(H,\S)$, we have $u(n)w=0$ for $n$
sufficiently large. Then define $V(H,\S)$ to be the quotient
$T(H\otimes \C[t,t^{-1}])$-module modulo the relation (\ref{e2.18})
for $u,v\in H,\; w\in \widetilde{V}(H,\S)$. Denote by ${\bf 1}$ the
image of $1+J_{+}$ in $V(H,\S)$. We see that $(V(H,\S),{\bf 1})$ is
a vacuum $(H,\S)$-module. Furthermore, we have (see \cite{li-qva2}):

\bp{pvacuum-existence} Let $H,\;\<\cdot,\cdot\>$ be given as before
and let $\S(x)$ be a linear map from $H\otimes H$ to $H\otimes
H\otimes \C[[x]]$. Then $(V(H,\S),{\bf 1})$ is a universal vacuum
$(H,\S)$-module. \ep

\br{rzfalgebra} {\em Let $H,\<\cdot,\cdot\>$ and $\S(x)$ be given as
in Proposition \ref{pvacuum-existence}. One can associate an
associative algebra $A(H,\S)$ as follows. First consider the tensor
algebra $T(H\otimes \C[t,t^{-1}])$ over the space $H\otimes
\C[t,t^{-1}]$. Then equip this algebra with the $\Z$-grading
uniquely determined by $\deg (u\otimes t^{n})=n$ for $u\in H,\; n\in
\Z$. Let $\overline{A(H,\S)}$ be the subspace of $\prod_{n\in
\Z}T(H\otimes \C[t,t^{-1}])_{n}$, consisting of those $f$ such that
$f(n)=0$ for $n$ sufficiently negative. Notice that
$\overline{A(H,\S)}$ is an associative algebra with $T(H\otimes
\C[t,t^{-1}])$ as a subalgebra. Then define $A(H,\S)$ to be the
quotient algebra of $\overline{A(H,\S)}$ modulo the relation
(\ref{e2.18}).} \er

The study on the weak quantum vertex algebras $V(H,\S)$ was
continued in \cite{kl} and it was proved that for certain special
cases, $V(H,\S)$ are irreducible quantum vertex algebras. Among
other results, the following were obtained in \cite{kl}:

\bt{tQqva} Let ${\bf Q}=(q_{ij}(x))$ be an $l\times l$ matrix in
$\C[[x]]$ such that $q_{ij}(x)q_{ji}(-x)=1$ for $1\le i,j\le l$ and
let $V$ be any nonzero nonlocal vertex algebra with a set of
generators $u^{(i)},v^{(i)}$ $(1\le i,j\le l)$ such that
\begin{eqnarray*}
&&Y(u^{(i)},x_{1})Y(u^{(j)},x_{2})=q_{ij}(x_{2}-x_{1})Y(u^{(j)},x_{2})Y(u^{(i)},x_{1}),\\
&&Y(v^{(i)},x_{1})Y(v^{(j)},x_{2})=q_{ij}(x_{2}-x_{1})Y(v^{(j)},x_{2})Y(v^{(i)},x_{1}),\\
&&Y(u^{(i)},x_{1})Y(v^{(j)},x_{2})-q_{ji}(x_{2}-x_{1})Y(v^{(j)},x_{2})Y(u^{(i)},x_{1})=
\delta_{ij}x_{2}^{-1}\delta\left(\frac{x_{1}}{x_{2}}\right).
\end{eqnarray*}
Then $V$ is an irreducible quantum vertex algebra. Furthermore, all
such quantum vertex algebras are isomorphic to each other. \et

\bp{pQqva-existence} Let ${\bf Q}=(q_{ij}(x))$ be an $l\times l$
matrix in $\C[[x]]$ as in Theorem \ref{tQqva}. Assume that there
exist $f_{ij}(x)\in \C[[x]]$ with $f_{ij}(0)\ne 0$ such that
$q_{ij}(x)=f_{ij}(x)f_{ji}(-x)^{-1}$ for $1\le i,j\le l$. Then
 there exists a nonzero irreducible quantum vertex
algebra as described in Theorem \ref{tQqva}. \ep

\bex{example-last} {\em Let $\lambda$ be any nonzero complex number.
In view of Theorem \ref{tQqva} and Proposition
\ref{pQqva-existence}, up to isomorphism there exists a unique
irreducible quantum vertex algebra $V[\lambda]$ with a set of
generators $u$ and $v$ such that
\begin{eqnarray*}
&&Y(u,x_{1})Y(u,x_{2})
=\iota_{x_{2},x_{1}}\left(\frac{\lambda+x_{1}-x_{2}}{\lambda-x_{1}+x_{2}}\right)
Y(u,x_{2})Y(u,x_{1}),\\
&&Y(v,x_{1})Y(v,x_{2})=\iota_{x_{2},x_{1}}
\left(\frac{\lambda+x_{1}-x_{2}}{\lambda-x_{1}+x_{2}}\right)Y(v,x_{2})Y(v,x_{1}),\\
&&Y(u,x_{1})Y(v,x_{2})-\iota_{x_{2},x_{1}}
\left(\frac{\lambda-x_{1}+x_{2}}{\lambda+x_{1}-x_{2}}\right)Y(v,x_{2})Y(u,x_{1})
=x_{2}^{-1}\delta\left(\frac{x_{1}}{x_{2}}\right).
\end{eqnarray*}
This can be viewed as a deformed $\beta\gamma$-system, where this
reduces to the standard $\beta\gamma$-system when $\lambda$ goes to
infinity.} \eex

Next, we discuss the example of quantum vertex algebras, constructed
by using a version of double Yangian $DY_{\hbar}(sl_{2})$. Let
$T(sl_{2}\otimes \C[t,t^{-1}])$ denote the tensor algebra over the
vector space $sl_{2}\otimes \C[t,t^{-1}]$. {}{\em From now on we
shall simply use ${\mathcal{T}}$ for this algebra.} We equip
${\mathcal{T}}$ with the $\Z$-grading which is uniquely defined by
$$\deg (u\otimes t^{n})=n\ \ \
\mbox{ for }u\in sl_{2},\; n\in \Z,$$ making ${\mathcal{T}}$ a
$\Z$-graded algebra ${\mathcal{T}}=\oplus_{n\in
\Z}{\mathcal{T}}_{n}$. Let $\overline{\mathcal{T}}$ denote the
subspace of $\prod_{n\in \Z}{\mathcal{T}}_{n}$, consisting of those
$f$ such that $f(n)=0$ for $n$ sufficiently negative. Then
$\overline{\mathcal{T}}$ is naturally an associative algebra with
$\mathcal{T}$ as a subalgebra.

For $u\in sl_{2}$, set
$$u(x)=\sum_{n\in \Z}u(n)x^{-n-1},$$
where $u(n)=u\otimes t^{n}$. We also write $sl_{2}(n)=sl_{2}\otimes
t^{n}$ for $n\in \Z$. Let $e,f,h$ be the standard Chevalley
generators of $sl_{2}$.

\bd{dyangian-double-algebra} {\em Let $q$ be a nonzero complex
number. Define $DY_{q}(sl_{2})$ to be the quotient algebra of
$\overline{\mathcal{T}}$ modulo the following relations:
\begin{eqnarray*}
& &e(x_{1})e(x_{2}) =\frac{q+x_{1}-x_{2}}{-q+x_{1}-x_{2}}
e(x_{2})e(x_{1}),\\
& &f(x_{1})f(x_{2}) =\frac{-q+x_{1}-x_{2}}{q+x_{1}-x_{2}}
f(x_{2})f(x_{1}),\\
& &[e(x_{1}),f(x_{2})]
=x_{1}^{-1}\delta\left(\frac{x_{2}}{x_{1}}\right)h(x_{2}),\\
& &h(x_{1})e(x_{2}) =\frac{q+x_{1}-x_{2}}{-q+x_{1}-x_{2}}
e(x_{2})h(x_{1}),\\
& &h(x_{1})h(x_{2})=h(x_{2})h(x_{1}),\\
& &h(x_{1})f(x_{2}) =\frac{-q+x_{1}-x_{2}}{q+x_{1}-x_{2}}
f(x_{2})h(x_{1}),
\end{eqnarray*}
where it is understood that
$$(\pm q+x_{1}-x_{2})^{-1}=\sum_{i\in \N}(\pm q)^{-i-1}(x_{2}-x_{1})^{i}
\in \C[[x_{1},x_{2}]].$$} \ed

As a convention, we define a $DY_{q}(sl_{2})$-{\em module} to be a
${\mathcal{T}}$-module $W$ such that for any $w\in W$,
$$sl_{2}(n)w=0 \ \ \ \ \mbox{ for $n$ sufficiently large}$$
and such that all the defining relations of $DY_{q}(sl_{2})$ in
Definition \ref{dyangian-double-algebra} hold. A vector $w_{0}$ of a
$DY_{q}(sl_{2})$-module $W$ is called a {\em vacuum vector} if
$$u(n)w_{0}=0\ \ \ \mbox{ for }u\in sl_{2},\; n\ge 0.$$
A {\em vacuum $DY_{q}(sl_{2})$-module} is a module $W$ equipped with
a vacuum vector $w_{0}$ that generates $W$. We use the pair
$(W,w_{0})$ to denote this vacuum module.

Set
$${\mathcal{T}}_{+}=\sum_{n\ge 1}{\mathcal{T}}_{n}\;\;\;\mbox{ and }\
J={\mathcal{T}}({\mathcal{T}}_{+}+sl_{2}(0)).$$ With $J$ a left
ideal of ${\mathcal{T}}$, ${\mathcal{T}}/J$ is a (left)
${\mathcal{T}}$-module and for any $v\in {\mathcal{T}}$,
$sl_{2}(n)(v+J)=0$ for $n$ sufficiently large. Then we define
$V_{q}$ to be the quotient ${\mathcal{T}}$-module of
${\mathcal{T}}/J$, modulo all the defining relations of
$DY_{q}(sl_{2})$. Denote by ${\bf 1}$ the image of $1+J$ in $V_{q}$.
{}From the construction, $(V_{q},{\bf 1})$ is a vacuum
$DY_{q}(sl_{2})$-module. It follows (\cite{li-infinity}, Lemma 2.5)
that $(V_{q},{\bf 1})$ is universal in the obvious sense. We have
the following result (see \cite{li-infinity}):

\bt{tmain-dy-0} For any nonzero complex number $q$, there exists a
(unique) weak quantum vertex algebra structure on $V_{q}$ with ${\bf
1}$ as the vacuum vector such that
$$Y(e(-1){\bf 1},x)=e(x),\ \ Y(f(-1){\bf 1},x)=f(x), \ \ \ Y(h(-1){\bf 1},x)=h(x).$$
Furthermore, $V_{q}$ is a non-degenerate quantum vertex algebra
 and for any $DY_{q}(sl_{2})$-module $W$, there exists
one and only one $V_{q}$-module structure $Y_{W}$ on $W$ such that
$$Y_{W}(e(-1){\bf 1},x)=e(x),\ \ Y_{W}(f(-1){\bf 1},x)=f(x),
\ \ Y_{W}(h(-1){\bf 1},x)=h(x).$$ \et

\section{Modules-at-infinity for quantum vertex algebras}
In Section 4, we have discussed a version  of double Yangian
$DY_{\hbar}(sl_{2})$ and its connection with quantum vertex
algebras. Here we discuss another version which is closer to the
original double Yangian (see \cite{kt}, \cite{kh}). {}From the
definition, modules for $DY_{q}(sl_{2})$ are of highest weight type.
In fact, the original double Yangians were defined by using a
completion corresponding to modules of lowest weight type. The fact
that quantum vertex algebras and their modules are of highest weight
type motivated us to introduce the algebra $DY_{q}(sl_{2})$ in
\cite{li-infinity}. On the other hand, the attempt to relate the
original double Yangian with quantum vertex algebras led us to a new
theory of modules-at-infinity for quantum vertex algebras and to a
conceptual construction of (weak) quantum vertex algebras and their
modules-at-infinity.

As in Section 4, let $T$ be the tensor algebra over the space
$sl_{2}\otimes \C[t,t^{-1}]$, equipped with the $\Z$-grading
uniquely determined by $\deg (sl_{2}\otimes t^{n})=n$ for $n\in \Z$.
 Let $\widetilde{T}$ be the subspace of $\prod_{n\in
\Z}T_{n}$, consisting of those $f$ such that $f(n)=0$ for $n$
sufficiently large.  Just as $\overline{T}$ in Section 4, the space
$\widetilde{T}$ is naturally an associative algebra containing $T$
as a subalgebra.

\bd{dyangian-infity} {\em Let $q$ be a nonzero complex number.
Define $DY_{q}^{\infty}(sl_{2})$ to be the quotient algebra of
$\widetilde{T}$ modulo the following relations:
\begin{eqnarray*}
& &e(x_{1})e(x_{2}) =\frac{x_{1}-x_{2}+q}{x_{1}-x_{2}-q}
e(x_{2})e(x_{1}),\\
& &f(x_{1})f(x_{2}) =\frac{x_{1}-x_{2}-q}{x_{1}-x_{2}+q}
f(x_{2})f(x_{1}),\\
& &[e(x_{1}),f(x_{2})]
=-x_{1}^{-1}\delta\left(\frac{x_{2}}{x_{1}}\right)h(x_{2}),\\
& &h(x_{1})e(x_{2}) =\frac{x_{1}-x_{2}+q}{x_{1}-x_{2}-q}
e(x_{2})h(x_{1}),\\
& &h(x_{1})h(x_{2})=h(x_{2})h(x_{1}),\\
& &h(x_{1})f(x_{2}) =\frac{x_{1}-x_{2}-q}{x_{1}-x_{2}+q}
f(x_{2})h(x_{1}),
\end{eqnarray*}
where it is understood that
$$(x_{1}-x_{2}\pm q)^{-1}=\sum_{i\in \N}(\pm q)^{i}(x_{1}-x_{2})^{-i-1}
\in \C[[x_{1}^{-1},x_{2}]]$$ (cf. Definition
\ref{dyangian-double-algebra}).} \ed

Let $W$ be a $T(sl_{2}\otimes \C[t,t^{-1}])$-module $W$ such that
for every $w\in W$,
\begin{eqnarray}
sl_{2}(n)w=0\ \ \ \mbox{ for $n$ sufficiently negative}
\end{eqnarray}
and such that all the defining relations for
$DY_{q}^{\infty}(sl_{2})$ hold. Then $W$ is naturally a
$DY_{q}^{\infty}(sl_{2})$-module. We call such a
$DY_{q}^{\infty}(sl_{2})$-module a {\em restricted} module. Note
that for any restricted $DY_{q}^{\infty}(sl_{2})$-module $W$, the
generating functions $e(x),f(x), h(x)$ are elements of the space
$\Hom (W,W((x^{-1})))$. In the following we discuss a conceptual
construction of weak quantum vertex algebras, by which $e(x),f(x),
h(x)$ generate a weak quantum vertex algebra.

Let $W$ be a general vector space over $\C$. Set
\begin{eqnarray}
\E^{o}(W)=\Hom (W,W((x^{-1})))\subset (\End W)[[x,x^{-1}]].
\end{eqnarray}
An (ordered) sequence $(\psi^{(1)}(x),\dots,\psi^{(r)}(x))$ in
$\E^{o}(W)$ is said to be {\em quasi-compatible} if there exists a
nonzero polynomial $p(x_{1},x_{2})$ such that
\begin{eqnarray}
\ \ \ \ \ \ \ \ \left(\prod_{1\le i<j\le r}p(x_{i},x_{j})\right)
\psi^{(1)}(x_{1})\cdots \psi^{(r)}(x_{r})\in \Hom
(W,W((x_{1}^{-1},\dots,x_{r}^{-1}))).
\end{eqnarray}
A subset $U$ of $\E^{o}(W)$ is said to be {\em quasi-compatible} if
every finite sequence in $U$ is quasi-compatible. In this
definition, if we assume that $p(x_{1},x_{2})$ is of the form
$(x_{1}-x_{2})^{k}$ with $k\in \N$, we then arrive at the notion of
{\em compatibility}.

Let $\iota_{x,\infty;x_{0},0}$ denote the field embedding
$$\iota_{x,\infty;x_{0},0}: \C(x,x_{0})\rightarrow
\C((x^{-1}))((x_{0})),$$ that preserves each polynomial. Let
$a(x),b(x)\in \E^{o}(W)$ such that $(a(x),b(x))$ is
quasi-compatible. By definition there exists a nonzero polynomial
$p(x_{1},x_{2})$ such that
\begin{eqnarray}\label{ecompatible-condition}
p(x_{1},x_{2})a(x_{1})b(x_{2})\in \Hom
(W,W((x_{1}^{-1},x_{2}^{-1}))).
\end{eqnarray}
As
\begin{eqnarray*}
& &\iota_{x,\infty;x_{0},0}\left(1/p(x_{0}+x,x)\right)\in
  \C((x^{-1}))((x_{0})),\\
& &\left( p(x_{1},x)a(x_{1})b(x)\right)|_{x_{1}=x+x_{0}} \in
(\Hom(W,W((x^{-1}))))[[x_{0}]],
\end{eqnarray*}
the product
$$\iota_{x,\infty;x_{0},0}\left(1/p(x_{0}+x,x)\right)
\left( p(x_{1},x)a(x_{1})b(x)\right)|_{x_{1}=x+x_{0}}$$ lies in
$(\Hom (W,W((x^{-1}))))((x_{0}))$.

\bd{doperation-same-order} {\em Let $(a(x),b(x))$ be a
quasi-compatible pair in $\E^{o}(W)$. We define
\begin{eqnarray*}
Y_{\E^{o}}(a(x),x_{0})b(x)=\iota_{x,\infty;x_{0},0}\left(1/p(x_{0}+x,x)\right)
\left( p(x_{1},x)a(x_{1})b(x)\right)|_{x_{1}=x+x_{0}},
\end{eqnarray*}
which lies in $\E^{o}(W)((x_{0}))$, where $p(x_{1},x_{2})$ is any
nonzero polynomial such that (\ref{ecompatible-condition}) holds.}
\ed

\bd{dmodule-at-infinity} {\em Let $V$ be a nonlocal vertex algebra.
A {\em $V$-module-at-infinity} is a vector space $W$, equipped with
a linear map $Y_{W}: V\rightarrow \Hom (W,W((x^{-1})))$, satisfying
the conditions that $Y_{W}({\bf 1},x)=1_{W}$ and that for $u,v\in
V$, there exists a nonnegative integer $k$ such that
$$(x_{1}-x_{2})^{k}Y_{W}(u,x_{1})Y_{W}(v,x_{2})\in \Hom
(W,W((x_{1}^{-1},x_{2}^{-1}))),$$ and
\begin{eqnarray}
\ \ \ \ \ \ x_{0}^{k}Y_{W}(Y(u,x_{0})v,x_{2})
=\left((x_{1}-x_{2})^{k}Y_{W}(u,x_{1})Y_{W}(v,x_{2})\right)|_{x_{1}=x_{2}+x_{0}}.
\end{eqnarray}}
\ed

In terms of this notion we have (see \cite{li-infinity}):

 \bt{tclosed} Let $V$ be a closed
compatible subspace of $\E^{o}(W)$, containing $1_{W}$. Then
$(V,Y_{\E^{o}},1_{W})$ carries the structure of a nonlocal vertex
algebra with $W$ as a (faithful) module-at-infinity where
$Y_{W}(\alpha(x),x_{0})=\alpha(x_{0})$. \et

\bt{tmaximal} Every maximal compatible subspace of $\E^{o}(W)$ is
closed and contains $1_{W}$. Furthermore, for any compatible subset
$U$, there exists a (unique) smallest closed compatible subspace
$\<U\>$ containing $U$ and $1_{W}$, and $(\<U\>,Y_{\E^{o}},1_{W})$
carries the structure of a nonlocal vertex algebra with $W$ as a
(faithful) module-at-infinity where
$Y_{W}(\alpha(x),x_{0})=\alpha(x_{0})$. \et

Recall that $\C(x)$ denotes the field of rational functions.  Let
$$\iota_{x,\infty}: \C(x)\rightarrow \C((x^{-1}))$$
be the unique field embedding that preserves each polynomial. To
obtain weak quantum vertex algebras we need analogs of $\S$-local
subsets.

\bd{ds-locality} {\em A subspace $U$ of $\E^{o}(W)$ is said to be
{\em $\S$-local} if for any $a(x),b(x)\in U$, there exists
$$\sum_{i=1}^{r}u^{(i)}(x)\otimes v^{(i)}(x)\otimes f_{i}(x)
\in U\otimes U\otimes \C(x)$$ such that
\begin{eqnarray*}\label{eslocality-relation}
(x_{1}-x_{2})^{k}a(x_{1})b(x_{2})=(x_{1}-x_{2})^{k}
\sum_{i=1}^{r}\iota_{x,\infty}(f_{i})(-x_{1}+x_{2})
u^{(i)}(x_{2})v^{(i)}(x_{1})
\end{eqnarray*}
for some nonnegative integer $k$.} \ed

In terms of this notion we have (see \cite{li-infinity}):

\bt{tqva-main} Let $W$ be a vector space and let $U$ be an
$\S$-local subset of $\E^{o}(W)$. Then $U$ is compatible and the
nonlocal vertex algebra $\<U\>$ generated by $U$ is a weak quantum
vertex algebra with $W$ as a module-at-infinity. \et

Recall from Section 4 the quantum vertex algebra $V_{q}$. Then we
have (see \cite{li-infinity}):

\bt{tdy-infinity} Let $q$ be any nonzero complex number and let $W$
be any restricted $DY_{q}^{\infty}(sl_{2})$-module. There exists one
and only one structure of a $V_{q}$-module-at-infinity on $W$ with
$$Y_{W}(e(-1){\bf 1},x)=e(x),\ \ Y_{W}(f(-1){\bf 1},x)=f(x),\ \
Y_{W}(h(-1){\bf 1},x)=h(x).$$ \et

\bibliographystyle{amsalpha}

\begin{thebibliography}{EFK}

\bibitem[A1]{ang1}
Iana I. Anguelova, Symmetric polynomials and $H_D$-quantum vertex
algebras, arXiv:math/0603145.

\bibitem[A2]{ang2}
Iana I. Anguelova, Super-bicharacter construction of quantum vertex
algebras, arXiv:0708.3708.

\bibitem[AB]{ab}
Iana I. Anguelova and Maarten J. Bergvelt, $H_{D}$-Quantum vertex
algebras and bicharacters, preprint; arXiv:0706.1528 [math.QA].

\bibitem [BK]{bk}
B. Bakalov and V. Kac, Field algebras, {\em Internat. Math. Res.
Notices} {\bf 3} (2003) 123-159.

\bibitem[B1]{b-va}
R. E. Borcherds, Vertex algebras, Kac-Moody algebras, and the
Monster, {\it Proc. Natl. Acad. Sci. USA} {\bf 83} (1986) 3068-3071.

\bibitem[B2]{b-gva}
R. E. Borcherds, Vertex algebras, in ``{\em Topological Field
Theory, Primitive Forms and Related Topics}'' (Kyoto, 1996), edited
by M. Kashiwara, A. Matsuo, K. Saito and I. Satake, Progress in
Math., Vol. 160, Birkh\"auser, Boston, 1998, 35-77.

\bibitem[B3]{b-qva}
R. E. Borcherds, Quantum vertex algebras, {\em Taniguchi Conference
on Mathematics Nara'98,} Adv. Stud. Pure Math., 31, Math. Soc.
Japan, Tokyo, 2001, 51-74.

\bibitem[Dr]{drinfeld}
V. G. Drinfeld, A new realization of Yangians and quantized affine
algebras, {\em Soviet Math. Dokl.} {\bf 36} (1988) 212-216.

\bibitem[EFK]{efk}
P. Etingof, I. Frenkel and A. Kirillov Jr., {\em Lectures on
Representation Theory and Knizhnik-Zamolodchikov Equations},
Mathematical Surveys and Monographs, Vol. 58, Amer. Math. Soc.,
Providence, 1998.

\bibitem[EK]{ek}
P. Etingof and D. Kazhdan, Quantization of Lie bialgebras, V, {\em
Selecta Mathematica, New Series}, {\bf 6} (2000) 105-130.

\bibitem[Fa]{fad}
L. Faddeev, Quantum completely integrable models in field theory,
{\em Soviet Sci. Rev., Ser. C: Math. Phys. Rev.} {\bf 1}, Hawood
Academic Publ., 1990, pp. 107-155.

\bibitem [FR]{efr}
E. Frenkel and N. Reshetikhin, Towards deformed chiral algebras, in:
{\em Quantum Group Symposium, Proc. of 1996 Goslar conference}, ed.
by H.-D. Doebner and V. K. Dobrev, Heron Press, Sofia, 1997, 27-42.

\bibitem [FJ]{fj}
I. B. Frenkel and N.-H. Jing, Vertex representations of quantum
affine algebras, {\em Proc. Natl. Acad. Sci. USA} {\bf 85} (1988)
9373-9377.

\bibitem[FLM]{flm}
I. Frenkel, J. Lepowsky and A. Meurman, {\it Vertex Operator
Algebras and the Monster}, Pure and Appl. Math., {\bf Vol. 134},
Academic Press, Boston, 1988.


\bibitem [K]{kac}
V. G. Kac, {\it Vertex Algebras for Beginners}, University Lecture
Series, Vol. 10, Amer. Math. Soc., 1997.

\bibitem [KL]{kl}
M. Karel and H.-S. Li, Some quantum vertex algebras of
Zamolodchikov-Faddeev type, {\em Commun. Contemp. Math.}, to appear.

\bibitem [Kh]{kh}
S. Khoroshkin, Central extension of the Yangian Double, arXiv:
q-alg/9602031.

\bibitem [KT]{kt}
S. Khoroshkin and V. Tolstoy, Yangian double, {\em Lett. Math.
Phys.} {\bf 36} (1996) 373-402.

\bibitem[Li1]{li-local}
H.-S. Li, Local systems of vertex operators, vertex superalgebras
and modules, {\em J. Pure Appl. Alg.} {\bf 109} (1996) 143-195;
hep-th/9406185.

\bibitem[Li2]{li-g1}
H.-S. Li, Axiomatic $G_{1}$-vertex algebras, {\em Commun. Contemp.
Math.} {\bf 5} (2003) 281-327.

\bibitem[Li3]{li-simple}
H.-S. Li, Simple vertex operator algebras are nondegenerate, {\em J.
Alg.} {\bf 267} (2003), 199-211.

\bibitem[Li4]{li-vpa}
H.-S. Li, Vertex algebras and vertex Poisson algebras, {\em Commun.
Contemporary Math.} {\bf 6} (2004), 61-110.

\bibitem[Li5]{li-new}
H.-S. Li, A new construction of vertex algebras and quasi modules
for vertex algebras, {\em Adv. Math.} {\bf 202} (2006), 232-286.

\bibitem[Li6]{li-qva1}
H.-S. Li, Nonlocal vertex akgebras generated by formal vertex
operators, {\em Sel. Math., New Ser.} {\bf 11} (2005) 349-397.

\bibitem[Li7]{li-qva2}
H.-S. Li, Constructing quantum vertex algebras, {\em International
Journal of Mathematics} {\bf 17} (2006), 441-476.

\bibitem[Li8]{li-infinity}
H.-S. Li, Modular-at-infinity for quantum vertex algebras, {\em
Commun. Math. Phys.} {\bf 282} (2008) 819-864.

\bibitem[Li9]{li-hqva}
H.-S. Li, $\hbar$-adic quantum vertex algebras and their modules,
arXiv:0812.3156 [math.QA].

\bibitem[ZZ]{zam}
A. B. Zamolodchikov and Al. B. Zamolodchikov, Factorized
$S$-matrices in two dimensionals as the exact solutions of certain
relativistic quantum field theory models, {\em Ann. of Physics} {\bf
120} (1979) 253-291.

\end{thebibliography}

\end{document}